\documentclass[twoside]{amsart}
\usepackage{amssymb,amsmath,amscd}
\usepackage{amsthm}
\usepackage[spanish]{babel}
\usepackage[latin1]{inputenc}  %acentos normales mpr
\usepackage{amsfonts}
\usepackage{enumerate}
\usepackage{amsbsy}
\usepackage{geometry}
\usepackage{graphicx}
\usepackage{subfig}
\usepackage{epsfig}
\usepackage[all]{xy}

\theoremstyle{plainm}
\newtheorem{teo}{Teorema}
\newtheorem{lem}[teo]{Lema}
\newtheorem{prop}[teo]{Proposici\'on}
\newtheorem{cor} [teo] {Corolario}

\theoremstyle{remarkm}

\theoremstyle{definitionm}
\newtheorem{dfn}[teo]{Definici\'on}

\newcommand{\bb}[1]{\mathbb{#1}}
\newcommand{\fk}[1]{\mathfrak{#1}}
\newcommand{\ml}[1]{\mathcal{#1}}

\begin{document}
	
	%%El siguiente pÂrrafo es el encabezado de la primera pÂgina del artÂculo, nosotros modificamos el nÂmero de pÂgina y el nÂmero de artÂculo
	%\setcounter{page}{9} \makeatletter
	%\def\cap@leftlogo{\vbox to\ht0{%
			%                  \fontspecialpage\footnotesize\hbox{\@themaname}%
			%                  \hbox{\vbox{ \hsize .5\hsize\noindent\strut%
					%                  \footnotesize\mixbaal@volnumber
					%                  \thepage\ --\ \pageref{artuno}
					%                  \strut}} \vfill}}
	%\makeatother
	%
	
	%%%%%%%%%%%%%%%%%%%%%%%%%%%%%%

	%Ejemplo
	\title[La curva de Fargues-Fontaine]{ La curva de Fargues--Fontaine: \\ Una motivaci\'on al estudio de la teor\'ia de representaciones de Galois $p$-\'adicas }
	\author[Jorge A. Robles Hdez. y J. Rogelio P\'erez B.]{Jorge Alberto Robles Hernandez\\
		 Jes\'us Rogelio P\'erez Buend\'ia}
	
	%%%%%%%%%%%%%%%%%%%%%%%%%%%%%%%
	\begin{abstract}
		Este art\'iculo proporciona una revisi\'on comprensiva sobre la curva de Fargues-Fontaine, una pieza central en la teor\'ia de Hodge $p$-\'adica, y su papel crucial en la clasificaci\'on de las representaciones de Galois $p$-\'adicas. Nos enfocamos en sintetizar los desarrollos fundamentales en torno a esta curva, subrayando c\'omo conecta conceptos avanzados de geometr\'ia aritm\'etica con la teor\'ia pr\'actica de representaciones. Analizamos en detalle los anillos de periodos de Fontaine ($B_{cris}, B_{st}, B_{dR}$), abordando sus propiedades algebraicas y aritm\'eticas esenciales, y c\'omo estos anillos contribuyen a la construcci\'on y definici\'on de la curva. Adem\'as, exploramos la teor\'ia de las representaciones de Galois $p$-\'adicas admisibles y discutimos c\'omo, una vez definida la curva, esta se relaciona con la teor\'ia de Harder-Narasimhan. 
	\end{abstract}
	
	\keywords{Curva de Fargues-Fontaine, Geometr\'ia Aritm\'etica, Representaciones de Galois $p$-\'adicas, Anillos de Peridos, Teoría de Hodge $p$-\'adica.}
	\subjclass{11F85, 11S15, 11S20, 11F80.}

	\maketitle
	
	\section{ Introducci\'on} 
	
	La geometr\'ia aritm\'etica ha experimentado un impulso significativo en los \'ultimos a\~nos gracias al uso de m\'etodos $p$-\'adicos. Uno de sus grandes avances ha sido el descubrimiento de la curva de Fargues--Fontaine que ha generado gran inter\'es debido a las diversas relaciones que tiene con distintas \'areas de las matem\'aticas~\cite{Fon3}.
	
	\vspace{.1cm}
	
	Fue descubierta en el a\~no 2009 por los matem\'aticos franceses Laurent Fargues y Jean Marc-Fontaine mediante una serie de correos electr\'onicos y en el transcurso de una semana de conferencias en Trieste, Italia. Pierre Colmez documenta este intercambio de mensajes y nos aporta el contexto mediante el cual se realiz\'o dicho descubrimiento~\cite{Far}. Fargues y Fontaine discut\'ian un art\'iculo de Berger~\cite{berger2008construction} donde aseguraba que el anillo de periodos $B_e$ era un anillo de Bezout. Fontaine, incr\'edulo de esta afirmaci\'on, indaga a\'un m\'as llegando a una conclusi\'on inesperada, el hecho que es de ideales principales. Fascinado por este anillo, Fontaine plantea a Colmez y Fargues una serie de preguntas cuyas respuestas desembocar\'ian en relaciones con la teor\'ia de filtraciones de Harder-Narasimhan en el contexto $p$-\'adico. Estas relaciones dar\'ian pi\'e al descubrimiento de la curva de Fargues--Fontaine revel\'andonos una conexi\'on importante, la relaci\'on de las fibras vectoriales de la curva con la teor\'ia de Harder-Narasimhan, que a su vez, tiene relaci\'on con las representaciones de Galois $p$-\'adicas y los llamados anillos de periodos de Fontaine.
	
	La teor\'ia de Hodge $p$-\'adica provee una forma de clasificar representaciones de Galois $p$-\'adicas de campos locales de caracter\'istica cero con campo residual de caracter\'istica $p$ (tambi\'en llamados campos de caracter\'istica mixta). Esta teor\'ia tiene sus inicios en los trabajos de Serre y Tate quienes estudiaban los m\'odulos de Tate sobre variedades abelianas y las representaciones de Hodge-Tate. Estas representaciones est\'an relacionadas a ciertas descomposiciones de cohomolog\'ias $p$-\'adicas, an\'alogas a la descomposici\'on de Hodge. Futuros avances en el \'area se inspiraron por las representaciones de Galois $p$-\'adicas que nacen de la cohomolog\'ia \'etale de variedades. Es en este contexto donde Jean-Marc Fontaine introduce los anillos de periods.
	
	% En el \'ambito de los m\'etodos $p$-\'adicos, es notable mencionar investigaciones sobre superficies K3 $p$-\'adicas llevadas a cabo por distintos grupos de investigaci\'on en M\'exico. Entre estas, se incluyen los trabajos que forman parte de la tesis doctoral del segundo autor \cite{perez2014crystalline}. Estas investigaciones destacan por su relevancia en la teor\'ia $p$-\'adica de Hodge y en el uso general de m\'etodos $p$-\'adicos. Como parte de nuestro compromiso con la educaci\'on en esta \'area, hemos impartido cursos de posgrado y dirigido investigaciones de tesis, incluyendo la actual del primer autor de este texto. Resaltamos, en particular, la escuela CIMPA ``Hodge Theory and P-adic Hodge Theory 2021'', organizada en CIMAT Guanajuato, M\'exico, donde el segundo autor imparti\'o un curso sobre representaciones de Galois $p$-\'adicas, disponible en l\'inea \cite{youtubeCIMPAHodgep}. Este curso, junto con otros recursos educativos disponibles en \cite{YoutubeCIMPASchoolPHodge}, constituyen una base s\'olida para aquellos interesados en profundizar en estos temas.
	
	La estructura del art\'iculo es el siguiente: En la secci\'on 2 presentamos notaci\'on e informaci\'on preliminar que se usar\'a a lo largo del texto. En la secci\'on 3 se plantea una introducci\'on a la idea de Fontaine sobre la clasificaci\'on de representaciones de Galois $p$-\'adicas para dar paso en la secci\'on 4 a la definici\'on de los anillos periodos y mencionar sus propiedades algebraicas y aritm\'eticas m\'as importantes. En la secci\'on 5 retomamos lo comentado en la secci\'on 3 con m\'as detalle para explicar de manera m\'as precisa el papel que toman los anillos de periodos en la clasificaci\'on de las representaciones de Galois $p$-\'adicas. En la secci\'on 6 presentamos la construcci\'on de la curva de Fargues-Fontaine usando como analog\'ia las ideas para la construcci\'on de la esfera de Riemann. Finalmente, en la secci\'on 7 describimos la conexi\'on entre las fibras vectoriales de la curva de Fargues-Fontaine con las representaciones de Galois $p$-\'adicas mediante el teorema de Harder-Narasimhan.
	
	\section{ Notaci\'on y Preliminares}
	
	Sea $p$ un primo fijo. El campo de los n\'umeros $p$-\'adicos (denotado como $\bb Q_p$) se puede definir como sigue: Para $a\in\bb Z$ sea $v_p(a)$ la m\'axima potencia de $p$ que divide al entero $a$. Extendiendo esto a $\bb Q$ definimos $v_p(\tfrac{a}{b})=v_p(a)-v_p(b)$. As\'i podemos definir una norma en $\bb Q$ como $$\left|\tfrac{a}{b}\right|_p = \tfrac{1}{p^{v_p\left(\tfrac{a}{b}\right)}}.$$

	Dicha norma induce una distancia en $\bb Q$ la cual no es completa (no toda sucesi\'on de Cauchy converge). A la completaci\'on de $\bb Q$ respecto a esta distancia es lo llamados el campo de los n\'umero $p$-\'adicos.
	
	Denotamos por $K$ a una extensi\'on finita de $\bb Q_p$, $\ml O_K$ su anillo de enteros, $m_K$ su \'unico ideal maximal y $k:=\ml O_K / m_K$ su campo residual. Denotamos por  $ K_0= K\cap \bb Q_p^{nr}$ a la m\'axima extensi\'on no ramificada~\cite{Ser} de $\bb Q_p$ dentro de $K$. Fijamos una cerradura algebraica $\overline{K}$ de $ K$ y sea $G_{K} = Gal(\overline{K}/K)$ su grupo de Galois absoluto. Notar que $\overline{K}$ es tambi\'en una cerradura algebraica de $\bb Q_p$ y por lo tanto no depende $K$.
	Sea $ K_0^{nr}$ la extensi\'on maximal no ramificada de $ K_0$ en $\overline{K}.$ An\'alogamente definimos a $ K^{nr}\subseteq \overline{K}.$ Como $ K_0$ es no ramificado sobre $\bb Q_p,$ $ K_0^{nr}$ es la extensi\'on m\'axima no ramificada de $\bb Q_p$ en $\overline{K}$ y por lo tanto tambi\'en es independiente de $K$.
	
	Denotemos por $\bb C_p$ a la completaci\'on $p-$\'adica de $\overline{K}$ y por $\ml O_{\bb C_p}$ a su anillo de enteros.  Fijamos $\pi\in m_K$ al par\'ametro uniformizador de $K$, es decir, un generador del ideal principal $m_K$. 
	
	Se dir\'a que un campo $K$ de caracter\'istica $p$ es perfecto si el morfismo de Frobenius es un isomorfismo.

	\section{ La idea de Fontaine.}
	
	Desde finales de la d\'ecada de 1970, J.M. Fontaine desarroll\'o un programa destinado a clasificar y describir las $\bb Q_p-$representaciones del grupo de Galois absoluto, $G_K$, de una extensi\'on finita de los racionales $p$-\'adicos $\mathbb Q_p$, i.e., los $\bb Q_p$ espacios vectoriales de dimensi\'on finita dotados de una acci\'on $\bb Q_p-$lineal continua de $G_K$~\cite{Bri}.
	
	La estrategia de Fontaine parte de la siguiente observaci\'on: si tenemos un anillo topol\'ogico $B$, dotado de una acci\'on continua de $G_K$ y estructuras adicionales estables bajo la acci\'on de $G_K,$ podemos asociar a cualquier representaci\'on de $G_K$ un invariante $D_B(V):=(B\otimes V)^{G_K}$ de $B\otimes V$.
	
	Entonces $D_B(V)$ es un $B^{G_K}$--m\'odulo equipado con estructuras adicionales heredadas de $B$, y que es a menudo m\'as f\'acil de describir que la representaci\'on $V$ de la que partimos. El anillo $B$ permite descomponer la sub--categor\'ia de $B-$representaciones admisibles (aquellas para los cuales \(B\otimes V\) es trivial, i.e., isomorfa a $B^{\dim V}$, como $G_K$- representaci\'on). Tales anillos topol\'ogicos $B$ son los llamados anillos de periodos de Fontaine. A continuaci\'on presentamos la construcci\'on de dichos anillos as\'i como de algunas de sus propiedades.
	
	\section{\hspace{.1cm} Teor\'ia de Hodge $p$-\'adica} \label{sec:style}
	
	La teor\'ia de Hodge $p$-\'adica, como lo describe \cite{serin2020}, puede verse desde dos puntos de vista: el aritm\'etico y el geom\'etrico.
	
	Desde el punto de vista aritm\'etico, es el estudio de las representaciones de Galois $p$-\'adicas, es decir, representaciones continuas $G_K\to Gl_n(\bb Q_p)$ donde $K$ es extensi\'on finita de $\bb Q_p$. Espec\'ificamente, esta teor\'ia busca construir un diccionario que relacione buenas categor\'ias de representaciones de $G_K$ con categor\'ias de objetos algebraicos semilineales. Un ejemplo de esto es el estudio de los m\'odulos de Tate de una curva el\'iptica sobre $K$ (que son representaciones de $G_K$) con buena reducci\'on junto con los llamados isocristales (es decir, $\bb Q_p$-espacios vectoriales de dimensi\'on finita equipados con un automofismo de Frobenius semilineal).
	
	Desde el punto de vista geom\'etrico, la teor\'ia de Hodge $p$-\'adica es el estudio de la geometr\'ia de una variedad (suave) $X$ sobre un campo $p$-\'adico $K$. En particular estamos interesados es varias teor\'ias de cohomolog\'ia relacionadas a $X$ como la cohomolog\'ia \'etale ($H_{et}^n(X_{\overline K}, \bb Q_p)$, la cohomolog\'ia de DeRham ($H_{dR}^n(X,K)$) y la cohomolog\'ia cristalina ($H_{cris}^n(X,K)$). Uno de los resultados mas relevantes de esta teor\'ia es derivado del caso cl\'asico en $\bb C$ sobre la descomposici\'on de Hodge de una variedad suave $Y$:
	$$  H^n (Y(\bb C,\bb Q)\otimes_{\bb Q}  \bb C \cong  \oplus_{i+j=n} H^i(Y,\Omega_Y^j). $$
	Tate observ\'o que exist\'ia una descomposici\'on an\'aloga para la cohomolog\'ia \'etale de una variedad abeliana sobre $K$ con buena reducci\'on, lo que lo llevo a conjeturar (y tiempo despu\'es Faltings lo demostrar\'ia) la descomposici\'on de Hodge-Tate \cite{faltings1988p} :
	
	$$  H^n_{et} (X_{\overline{K}},\bb Q_p)\otimes_{\bb Q_p}  \bb C_K \cong  \oplus_{i+j=n} H^i(X,\Omega_{X/K}^j) \otimes_{K}  \bb C_K(-j),  $$
	donde $X$ es una variedad suave sobre $K$ y es un isomorfismo de representaciones de Galois $p$-\'adicas $Rep_{\bb Q_p}(G_K)$.
	
	En este sentido, Fontaine realiza una serie de conjeturas, ahora teoremas, conocidos como \textgravedbl teoremas de comparaci\'on\textacutedbl en donde relaciona a las distintas cohomolog\'ias de la variedad, siendo necesario extender los coeficientes a los llamados anillos de periodos de Fontaine. 
	
	\begin{teo} \cite{faltings1988p}
		Sea $K$ una extensi\'on finita de $\bb Q_p$ y $X$ una variedad proyectiva lisa definida sobre $K$. Existen isomorfismos naturales:
		\begin{itemize}
			\item $H^n_{et} (X_{\overline{K}},\bb Q_p)\otimes_{\bb Q_p}  B_{cris} \cong  H^i_{cris}(X) \otimes_{K_0}   B_{cris}$.
			\item  $H^n_{et} (X_{\overline{K}},\bb Q_p)\otimes_{\bb Q_p}  B_{st} \cong  H^i_{HK}(X) \otimes_{K_0}   B_{st}$.
			\item $H^n_{et} (X_{\overline{K}},\bb Q_p)\otimes_{\bb Q_p}  B_{dR} \cong  H^i_{dR}(X) \otimes_{K}   B_{dR}$.
		\end{itemize}
		que conmutan con las acciones de $G_K$, Frobenius, $N$ (monodrom\'ia) y respectivas filtraciones. Aqu\'i $H^i_{HK}$ denota la cohomolog\'ia de Hyodo-Kato~\cite{hyodo1994semi}.
	\end{teo}
	
	En esta secci\'on nos centraremos en estudiar los anillos de periodos, en particular, mostraremos sus construcciones y las propiedades aritm\'eticas y algebraicas que nos servir\'an para la construcci\'on de la curva de Fargues-Fontaine. En el presente texto no se estudiar\'a la relaci\'on con las cohomolog\'ias como en el teorema 4.1. Para el lector interesado en este aspecto puede consultar \cite{Bri}.
	
	\subsection{ Anillos de Periodos}
	
	En geometr\'ia algebraica la palabra \textgravedbl periodo\textacutedbl se suele referir a un n\'umero complejo que puede ser expresado como integral de una funci\'on algebraica sobre un  dominio algebraico~\cite{kontsevich2001periods}. Uno de ellos es $2i\pi=\int_\gamma \frac{dt}{t},$ donde $\gamma$ es el c\'irculo unitario en el plano complejo. La teor\'ia de Hodge $p$-\'adica nos permite dise\~nar un an\'alogo $p$-\'adico de los periodos, trabajo que realiz\'o Fontaine creando anillos espec\'ificos~\cite{fontaine1994representations,Fon1,Fon2,Fon3}, los cuales adem\'as est\'an estrechamente relacionados con distintos tipos de cohomolog\'ias. Estos anillos forman parte fundamental en la construcci\'on de la curva de Fargues-Fontaine  y las representaciones de Galois $p$-\'adicas. 
	
	En esta secci\'on se dar\'a una introducci\'on a los anillos de periodos como en~\cite{Car}, en particular estudiaremos a los anillos $B_{cris}$, $B_{st}$ y $B_{dR}$, sus respectivas construcciones y algunas de sus propiedades algebraicas y anal\'iticas.
	
	\subsubsection{ El anillo  $B_{inf}$}
	
	Denotemos como $\phi$ el morfismo de Frobenius $x\mapsto x^p$ actuando en el cociente $\frac{\ml O_{\bb C_p} }{p\ml O_{\bb C_p} }$ y notemos que $\phi$ es un homomorfismo de anillos. Sea $R$ el l\'imite  del sistema proyectivo de anillos:
	
	$$\xymatrix{
		\frac{\ml O_{\bb C_p} }{p\ml O_{\bb C_p}} & \ar_{\phi} [l]  \frac{\ml O_{\bb C_p} }{p\ml O_{\bb C_p}}  & \ar_{\phi} [l] \cdots & \ar_{\phi} [l]  \frac{\ml O_{\bb C_p} }{p\ml O_{\bb C_p}}  \cdots
	},$$
	
	Espec\'ificamente, un elemento de $R$ es una sucesi\'on $(\zeta_n)_{n\geq 0}$ de elementos en $\tfrac{\ml O_{\bb C_p} }{p\ml O_{\bb C_p} }$, que satisfacen la propiedad de compatibilidad $\zeta_{n+1}^p =\zeta_n, \hspace{.2cm} \forall n\geq 0.$ $R$ es un anillo perfecto de caracter\'istica $p$.
	
	El anillo $R$ esta equipado con una valuaci\'on $v_b$ que se definir\'a a continuaci\'on: notemos que si $x\in \tfrac{\ml O_{\bb C_p} }{p\ml O_{\bb C_p} } \backslash \{0\},$ la valuaci\'on $p$--\'adica de $\widehat{x}\in \ml O_{\bb C_p}$ no depende del levantamiento $\widehat{x}$ de $x.$ La valuaci\'on $v_p:\ml O_{\bb C_p}$ induce una funci\'on bien definida en el cociente \(v_p:\tfrac{\ml O_{\bb C_p} }{p\ml O_{\bb C_p} } \to \bb Q \cup \{+\infty\}\) donde $v_p(0)=+\infty.$ As\'i, para $\zeta=(\zeta_n)_{n\geq 0}$ en $R$ definimos
	$$v_b(\zeta) := \lim_{n\to\infty} p^n v_p(\zeta_n).$$
	La condici\'on de compatibilidad $\zeta_{n+1}^p=\zeta_n$ Implica que la sucesi\'on $(p^nv_p(\zeta_n))_{n\geq0}$ es eventualmente constante y el l\'imite est\'a bien definido.
	\vspace{.2cm}
	
	A continuaci\'on se har\'a uso de una estructura algebraica conocida como el anillo de vectores de Witt  y cuya idea es la de construir extensiones no ramificadas.
	\begin{dfn}
		Definimos como $A_{inf}=W(R)$ donde $W(-)$ es el funtor de vectores de Witt y $B_{inf}^+= A_{inf}\left[\tfrac{1}{p}\right] = W(R)\left[\tfrac{1}{p}\right],$ la localizaci\'on en $A_{inf}$ en $p.$
	\end{dfn}
	
	Para $x\in R,$ definimos el representante de Techm\"{u}ller en $A_{inf}$ como  $[x]=(x,0,\cdots,0,\cdots)$. Dado que $R$ es de valuaci\'on discreta, todo elemento de $A_{inf}$ puede ser escrito de manera \'unica como $$\sum_{i\geq 0} [\zeta_i]p^i, \quad \zeta_i\in R.$$ 
	
	De manera similar, todo elemento de $B_{inf}^+$ se puede escribir de manera \'unica de la forma $$\sum_{i\geq i_0}[\zeta_i]p^i, \quad \zeta_i\in R,$$ donde $i_0$ puede ser negativo y depende de $x.$
	
	Adem\'as, $B_{inf}^+$ cuenta con estructuras adicionales:
	\begin{itemize}
		\item $B_{inf}^+$ tiene una acci\'on de Frobenius $\varphi$ dado por:
		$$\varphi\bigg(\sum_{i=i_0}^\infty [\zeta_i]p^i \bigg) =\sum_{i=i_0}^\infty [\zeta_i^p]p^i,\hspace{.5cm} \zeta_i\in R, \hspace{.2cm} i_0\in\bb Z.$$
		
		\item $B_{inf}^+$ est\'a equipado con una acci\'on de $G_K$ dado por:
		$$g\bigg(\sum_{i=i_0}^\infty [\zeta_i]p^i \bigg) =\sum_{i=i_0}^\infty [g\zeta_i]p^i,\hspace{.5cm} \zeta_i\in R, \hspace{.2cm} i_0\in\bb Z, \hspace{.2cm} \forall g\in G_K.$$
	\end{itemize}

	Fijamos $\epsilon_1$ una ra\'iz primitiva $p$--\'esima de la unidad en $\ml O_{\overline{ K}}.$ Escojamos a $\epsilon_2$ como ra\'iz primitiva $p$--\'esima de $\epsilon_1.$ Entonces $\epsilon_2$ es ra\'iz primitiva $p^2$--\'esima de la unidad. Repitiendo este proceso, construimos elementos $\epsilon_3,\epsilon_4,...\in\ml O_{\overline{ K}}$ tales que $\epsilon_{n+1}^{p}=\epsilon_n, \hspace{.2cm} \forall n\geq 0.$ Sea $\overline{\epsilon}_n\in \frac{\ml O_{\overline{ K}}}{p\ml O_{\overline{ K}}}$ la clase de $\epsilon_n.$ Por la propiedad de compatibilidad, $\underline{\epsilon} = (1,\overline{\epsilon}_1,\overline{\epsilon}_2,...)\in R.$
	
	\begin{dfn}
		Para $\zeta=(c_0,c_1,...)\in R,$ definimos $\zeta^{\sharp}=\lim_{n\to\infty}\widehat{c_n}^{p^n},$ donde $\widehat{c_n}$ es un levantamiento de $c_n$ en $\ml O_{\bb C_p}$. 
	\end{dfn}
	
	La funci\'on $\sharp:R\to\ml O_{\bb C_p}, \hspace{.2cm} \zeta \mapsto \zeta^{\sharp} $ es inyectiva y multiplicativa. Por las propiedades de los vectores de Witt, la funci\'on \textgravedbl sharp\textacutedbl se extiende a un homomorfismo inyectivo de $\widehat{K}_0^{ur}-$\'algebras $\theta: B_{inf}^{+} \to\bb C_p$ que conmuta con la acci\'on de $G_K.$ Est\'a dado por
	$$\sum_{i=i_0}^\infty[\zeta_i]p^i\longmapsto\sum_{i=i_0}^\infty \zeta_i^{\sharp}p^i,\hspace{.2cm} i_0\in\bb Z, \hspace{.2cm} \zeta_i\in R.$$
	
	La siguiente proposici\'on nos dice que el kernel de $\theta$ es principal y nos muestra expl\'icitamente a un generador del mismo.
	
	\begin{prop}[\cite{Car}, Lema:3.1.5]
		Sea $z\in A_{inf}$ elemento tal que $\theta(z)=0$ y $v_b(z\pmod p)=1.$ Entonces $z$ genera $A_{inf}\cap \ker\theta$ como ideal de $A_{inf}.$ En particular, el elemento
		$$\omega = \dfrac{[\underline{\epsilon}]-1}{[\underline{\epsilon}^{\frac{1}{p}}]-1} = [\underline{\epsilon}^{\frac{1}{p}}]+[\underline{\epsilon}^{\frac{1}{p}}]^2+ \cdots + [\underline{\epsilon}^{\frac{1}{p}}]^{p-1}, $$
		satisface las condiciones.
	\end{prop}
	
	\subsubsection{\hspace{.6cm} El anillo $B_{cris}$}
	
	Para definir el anillo $B_{cris}$ necesitamos tener la noci\'on de elementos de la forma $\tfrac{x^n}{n!}$ y para ello daremos un breve pre\'ambulo a las llamadas potencias divididas. Se considerar\'a $0!=1.$
	
	\begin{dfn}
		Sea $A$ un anillo (conmutativo con unidad), $I$ un ideal $A.$ Una colecci\'on de aplicaciones $\gamma_n:I \to I, \hspace{.2cm} n\geq0$ y $\gamma_0:A\to I$ definido como $\gamma_0(x)=1$, es llamado una estructura de potencias divididas en $I$ si para todos $ n\geq0, m>0, x,y\in I$ y $a\in A$ se tiene:
		\begin{itemize}
			\item [1)] $\gamma_1(x)=x$
			
			\item [2)] $\gamma_n(x)\gamma_m(x)= \dfrac{(n+m)!}{n!m!}\gamma_{n+m}(x).$
			
			\item [3)] $\gamma_n(ax)=a^n\gamma_n(x).$
			
			\item [4)] $\gamma_n(x+y)=\sum_{i=0}^n \gamma_i(x)\gamma_{n-i}(y).$
			
			\item [5)] $\gamma_n(\gamma_m(x)) = \dfrac{(nm)!}{n!(m!)^n} \gamma_{nm}(x).$
		\end{itemize}
	\end{dfn}
	
	Notemos que $\tfrac{(n+m)!}{n!m!} = \binom{n+m}{n} \in \bb Z$ y adem\'as, $\tfrac{(nm)!}{n!(m!)^n}\in \bb Z$ ya que cuenta el n\'umero de maneras de dividir un grupo de $nm$ objetos en $n$ grupos de $m$ elementos.
	
	\begin{lem}
		Sea $A$ un anillo, $I$ un ideal de $A$. Si $\gamma$ es una estructura de potencias divididas en $I$, entonces $n!\gamma_n(x)=x^n$ para todo $n\geq1$ y $x\in I$.
	\end{lem}
	
	Este lema se puede demostrar por inducci\'on. As\'i, en el anillo $A$ tendr\'iamos la noci\'on de dividir entre $n!$, cuesti\'on que no siempre es posible si $char(A)>0$.
	
	\vspace{.1cm}
	
	Volviendo a los anillos de periodos, dado $x\in A_{inf}=W(R)$ denotamos como $A_{inf}\langle x\rangle$ la sub $A_{inf}$--\'algebra de $B_{inf}^+$ generada por las potencias divididas.

	\begin{dfn}
		Definimos $A_{cris}$ como la completaci\'on $p$--\'adica de $A_{inf}\langle z\rangle,$ donde $z$ es un generador de $A_{inf}\cap \ker\theta.$ Denotamos por $B_{cris}^{+}= A_{cris}\left[\frac{1}{p}\right]$ a la localizaci\'on en $p$ de $A_{cris}$. 
	\end{dfn}
	Como $\omega$ es un generador de $A_{inf}\cap \ker\theta$, se tiene que $A_{cris}= \widehat{ A}_{inf}  \langle \omega\rangle$, donde $\widehat{\phantom{A}}$ denota la completaci\'on $p$--\'adica. 
	Otro generador de $A_{inf}\cap \ker \theta$ tambi\'en es el elemento $[p^b]-p$.
	
	El siguiente lema nos da una condici\'on sobre cu\'ando cualesquiera dos elementos en $A_{inf}$ producen la misma sub $A_{inf}$-\'algebra generada por las potencias divididas.
	
	\begin{lem}
		Si $x,y\in A_{inf},$ $x\equiv y \pmod {pA_{inf}}$ entonces $A_{inf}\langle x\rangle = A_{inf} \langle y\rangle. $
	\end{lem} 
	
	Usando el lema, se tendr\'ia que $A_{cris}=\widehat{A}_{inf}\langle [p^b]\rangle$.
	
	Como $A_{cris}$ est\'a definida como completaci\'on $p$--\'adica es natural equipar $A_{cris}$ y $B_{cris}^+$ con la topolog\'ia $p$--\'adica. Con esta topolog\'ia, la inclusi\'on $A_{inf}\hookrightarrow A_{cris}$ es continua al igual que la inclusi\'on $B_{inf}^+\hookrightarrow B_{cris}^+$.
	
	Frobenius se extiende can\'onicamente a un endomorfismo $A_{cris}\to A_{cris},$ ya que $A_{inf}\langle[p^b]\rangle$ es estable bajo Frobenius. En efecto,
	$$\varphi\bigg( \dfrac{[p^b]^n}{n!}\bigg) = [p^b]^{np-n}\varphi\bigg(\dfrac{[p^b]^n}{n!}\bigg). $$
	
	Invirtiendo $p$ se obtiene Frobenius para $B_{cris}^+.$ An\'alogamente, $G_K$ extiende la acci\'on a $B_{cris}^+.$ 
	
	$A_{cris}$ contiene un periodo para el caracter ciclot\'omico, i.e., un elemento en el cu\'al Galois act\'ua por multiplicaci\'on por $\chi$, el caracter ciclot\'omico. Este elemento es
	$$t=log([\underline{\epsilon}]) =\sum_{i=1}^\infty(-1)^{i-1} \dfrac{([\underline{\epsilon}]-1  )^i}{i}.$$
	
	Frobenius act\'ua como $[\underline{\epsilon}]\mapsto[\underline{\epsilon}]^p$ y $G_K$ act\'ua como $g[\underline{\epsilon}]=[\underline{\epsilon}]^{\chi(g)}, \hspace{.2cm}g\in G_K.$ Tomando logaritmos, $\varphi(t)=pt$ y $gt=\chi(g)t, \hspace{.2cm}\forall g\in G_K.$ As\'i, $t$ es un periodo del caracter ciclot\'omico.
	
	\subsubsection{ El anillo $B_{dR}$}
	
	\begin{dfn}
		Definimos $B_{dR}^+$ como la completaci\'on de $B_{inf}^+$ respecto a la topolog\'ia $(\ker \theta)$--\'adica, es decir
		$$B_{dR}^+=\varprojlim_m \dfrac{B^+_{inf}}{(\ker\theta)^m}.$$
	\end{dfn}
	
	Como $B^+_{dR}$ est\'a definido como una completaci\'on, la topolog\'ia natural en este anillo es la topolog\'ia $(\ker\theta)-$\'adica. Una sucesi\'on de elementos $(x_n)_{n\geq 0}$ de elementos en $B^+_{dR}$ converge a $x\in B^+_{dR}$ si y s\'olo si para todo $m,$ la sucesi\'on $\{x_n\pmod{Fil^mB^+_{dR}}\}$ es eventualmente constante.
	
	\begin{dfn}
		Definimos el anillo $B_{dR}=B^+_{dR}[\frac{1}{t}].$
	\end{dfn}
	
	Como $B^+_{dR}$ es un anillo de valuaci\'on discreta con uniformizador $t$, se tiene que $B_{dR}$ es el campo de fracciones de $B^+_{dR}$, es decir, este anillo es de hecho un campo.
	
	La filtraci\'on de De Rham se extiende a $B_{dR}:$ $Fil^mB_{dR}=t^mB^+_{dR},\hspace{.2cm} m\in\bb Z$. Para m\'as informaci\'on se recomienda consultar al lector \cite{Fon1}, \cite{Fon2}, \cite{Car}. 
	
	\section{Representaciones de Galois $p$-\'adicas} 
	
	Una vez que hemos constru\'ido los anillos de periodos, regresamos a la idea de Fontaine para clasificar las representaciones de Galois $p$-\'adicas. En esta secci\'on se hablar\'a primero de la $B$-admisibilidad de una representaci\'on para luego dar paso a un ejemplo concreto: como el anillo de periodo $B_{dR}$ nos dar\'a m\'as informaci\'on sobre la representaci\'on.
	La bibliograf\'ia que se seguir\'a en esta secci\'on es \cite{Bri}.
	
	\subsection{Representaciones admisibles} 
	
	Sea $F$ un campo y $G$ un grupo. Sea $B$ una $F$--\'algebra, dominio equipado con una $G$--acci\'on (como $F$--\'algebra) y asuma que la sub $F$--\'algebra $E=B^G$ es un campo.
	
	No se imponen estructuras topol\'ogicas en $B,F$ o $G.$ El objetivo en usar a $B$ para construir un funtor de representaciones de $G$ $F$--lineales de dimensi\'on finita a $E$--espacios vectoriales de dimensi\'on finita equipados con estructuras adicionales que dependen de $B.$
	
	Sea $C=Frac(B)$ y $G$ act\'ua en $C$ de manera natural, i.e., $g(\frac{a}{b}) = g(a) g (b^{-1}).$
	
	\begin{dfn}
		Decimos que $B$ es $(F,G)$--regular si $C^G=B^G(=E)$ y si para todo $b\in B\backslash \{0\}$ cuyo espacio generado $F$--lineal $Fb$ es $G$--estable se tiene que $b$ es unidad en $B.$
	\end{dfn}
	
	Notemos que si $B$ es un campo, entonces es $(F,G)$--regular.
	
	\begin{dfn}
		Si $B$ es un dominio $(F,G)$--regular y $E$ den\'ota el campo $C^G=B^G,$ entonces para cualquier objeto $V\in Rep_F(G)$ de $G$ representaciones  $F$--lineales de dimensi\'on finita definimos $$D_B(V):=(B\otimes_F V)^{G},$$ es decir, $D_B(V)$ es un $E$--espacio vectorial equipado con un mapeo can\'onicoo
		$$\alpha_V: B\otimes_E D_B(V)\to  B\otimes_E (B\otimes_F V) = (B\otimes_E B)\otimes_F V \to B \otimes_F V.$$
	\end{dfn}
	
	\begin{dfn}
		Si se tiene la igualdad $\dim_E D_B(V)=\dim_F (V)$ se dir\'a que $V$ es una representaci\'on $B$--admisible.
	\end{dfn}
	
	El siguiente teorema nos indica que el $E$-espacio vectorial $D_B(V)$ es de hecho de dimensi\'on finita y que al restringirnos a las representaciones $B$-admisibles, se tiene un funtor exacto y fiel.
	
	\begin{teo}
		Fijamos $V$ como antes.
		\begin{itemize}
			\item [i)] El mapeo $\alpha_v$ es siempre inyectivo y $\dim_E D_B(V)\leq\dim_F V.$ La igualdad se da y si s\'olo si $\alpha_V$ es un isomorfismo.
			
			\item [ii)] Sea $Rep_F^B(G)\subseteq Rep_F(G)$ la subcategor\'ia de representaciones $B$--admisibles. El funtor contravariante $D_B:Rep_F^B(G)\to Vec_E$ es exacto y fiel.
		\end{itemize}
	\end{teo}
	
	En particular, considerando a $K$ un campo $p$--\'adico, i.e., extensi\'on finita de $\bb Q_p,$ $F=\bb Q_p$ y $G=Gal(\overline{K},K)$ tenemos que:
	\begin{itemize}
		\item Si $B=B_{dR}$, se dice que la representaci\'on es de DeRham.
		
		\item Si $B=B_{cris}$, se dice que la representaci\'on es cristalina.
		
		\item Si $B=B_{st}$, se dice que la representaci\'on es semiestable.
	\end{itemize}

	\subsection{Representaciones de deRham}
	
	Analizaremos un poco m\'as a detalle las representaciones de DeRham. Como $B_{dR}$ es $(\bb Q_p,G_K)$--regular con $B_{dR}^{G_K}=K$, la formalizaci\'on general de representaciones admisibles provee de una buena clase de representaciones $p$-\'adicas, los que son $B_{dR}$-admisibles.
	
	\begin{dfn}
		Definimos el funtor covariante $D_{dR}:Rep_{\bb Q_p}(G_K)\to Vec_K$ a la categor\'ia de $K$--espacios vectoriales de dimensi\'on finita como $D_{dR}(V)=(B_{dR}\otimes_{\bb Q_p}V)^{G_K}.$
		En el caso que $dim_{B_{dR}}(V)=dim_{\bb Q_p}(V)$ diremos que $V$ es una representaci\'on de deRham.
	\end{dfn}
	
	El codominio del funtor anterior tiene estructuras $K$-lineales adicionales (que vienen de la estructura adicional de la $K$-\'algebra $B_{dR}$), espec\'ificamente una filtraci\'on $K$-lineal que surge de la filtraci\'on $K$--lineal en el campo de fracciones $B_{dR}$ del anillo de valuaci\'on discreto completo $B_{dR}^+$ sobre $K.$
	
	Para $V\in Rep_{\bb Q_p}(G_K),$ el $K$-espacio vectorial $D_{dR}(V)=(B_{dR}\otimes_{\bb Q_p} V)^{G_K}\in Vec_K$ tiene estructura natural de objeto en $Fil_K:$ Como $B_{dR}$ tiene una filtraci\'on $K$--lineal, $G_K$ estable dado por $Fil^i(B_{dR})=t^{i}B_{dR}^+,$ obtenemos una filtraci\'on $K$--lineal, $G_K$--estable $\{ Fil^i(B_{dR})\otimes_{\bb Q_p}V \}$ en $B_{dR}\otimes_{\bb Q_p} V$ y este induce una filtraci\'on en $D_{dR}(V)$ de elementos $G_K$--invariantes, expl\'icitamente $$Fil^i(D_{dR}(V))= (t^iB_{dR}^+\otimes_{\bb Q_p}V)^{G_K}.$$
	
	\section{ Construcci\'on de la Curva de Fargues-Fontaine.}
	
	La curva de Fargues--Fontaine es un objeto de la teor\'ia de n\'umeros descubierta en 2009 por Laurent Fargues y Jean Marc--Fontaine, donde en su geometr\'ia codifica mucha informaci\'on sobre la aritm\'etica de los n\'umeros $p$--\'adicos. Se ha convertido r\'apidamente en un tema de investigaci\'on en la teor\'ia de Hodge $p$-\'adica y el programa de Langlands.
	
	En esta secci\'on motivamos la definici\'on de la curva usando una analog\'ia con la esfera de Riemann. Luego, construimos la curva desde dos puntos de vista distintos, uno desde el punto de vista del \'algebra conmutativa y otro desde los espacios de \textgravedbl Tilts y Untilts\textacutedbl \cite{Morr1, Col}.

	\subsection{La esfera de Riemann.}
	
	Para motivar la definici\'on de la curva, analizaremos una curva m\'as familiar, la esfera de Riemann $\bb P_{\bb C}^1$, a la cual le podemos asociar anillos $\bb C[z]\subseteq\bb C((\frac{1}{z})).$
	\begin{itemize}
		\item Primero, el anillo de funciones meromorfas sobre $\bb P_{\bb C}^1$ sin polos fuera del punto al infinito es el \'algebra de polinomios $\bb C[z],$ donde $z$ denota el par\'ametro local habitual en el origen.
		
		\item Mirando los desarrollos de Laurent en el punto al infinito de todas las funciones merom\'orficas en $\bb P_{\bb C}^1$ encontramos el anillo $\bb C((\frac{1}{z})),$ donde $\frac{1}{z}$ es un par\'ametro local para el punto al infinito.
	\end{itemize}
	
	    \begin{figure}[h]
	    	\includegraphics[scale=.4]{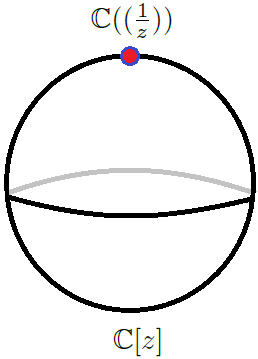}
	    	\caption{Esfera de Riemann.}
	    	\centering
	    \end{figure}
	
	Rec\'iprocamente, podemos reconstruir la esfera de Riemann a partir de los anillos $\bb C[z]\subseteq \bb C((\frac{1}{z})).$
	
	$$x\in\bb P_{\bb C}^1 \leftrightsquigarrow l_x\in \{  f\in \bb C[z] \hspace{.2cm}|\hspace{.2cm} \deg f\leq1 \} = \bb C \oplus \bb Cz.$$
	
	M\'as a\'un, si $S=\bigoplus_{k\geq0} \{ f\in \bb C[z] \hspace{.2cm}|\hspace{.2cm} \deg f\leq k \},$ entonces $Proj(S)=\bb P_{\bb C}^1.$\\
	
	La idea de la curva de Fargues--Fontaine es emular  la construcci\'on de la esfera de Riemann con los anillos $\bb C[z]\subseteq \bb C ((\frac{1}{z}))$; pero en los n\'umeros $p$--\'adicos y usando a los anillos $B_e\subseteq B_{dR},$ donde $B_e:=B_{cris}^{\varphi=1}.$
	
	\begin{figure}[h]
		\centering
		\includegraphics[scale=.15]{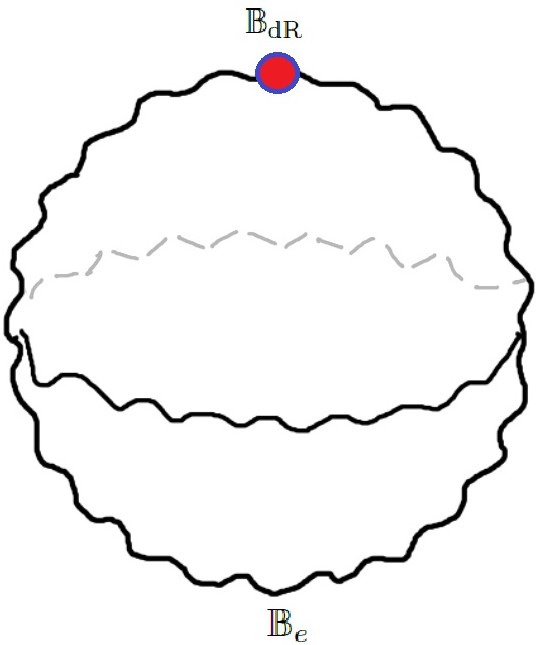}
		\caption{Representaci\'on art\'istica de la curva de Fargues--Fontaine \cite{Morr1}.}
	\end{figure}
	
	A esta curva la denotaremos como $X^{FF}.$ Concretamente, Fargues y Fontaine demostraron que \cite{Far}: 
	
	\begin{teo}[FF]
		Existe un esquema regular, noetheriano, conexo, separado de dimensi\'on uno, $X^{FF}$, sobre $\bb Q_p$;  tal que tiene un punto al infinito y, el anillo de funciones meromorfas sin polos fuera del infinito es el anillo $B_e.$
	\end{teo}
	
	En este caso, si $S=\bigoplus_{k\geq 0} \{f\in  B_e\hspace{.2cm}|\hspace{.2cm} \deg f\leq k \},$ entonces la curva de FF se define como  $X^{FF}=Proj(S).$ 
	Adem\'as, de manera an\'aloga al caso de la esfera de Riemann, tenemos que: 
	$$x\in X^{FF} \leftrightsquigarrow \{f\in B_e\hspace{.2cm}|\hspace{.2cm} \deg f\leq 1 \} = l_x.$$
	
	\subsection{Otro punto de vista: Tilts y Untilts.}
	
	A continuaci\'on presentamos otra forma de construir la curva de Fargues-Fontaine mediante lo llamados tilts y untilts. Esta secci\'on est\'a basado en la secci\'on 2.1 de \cite{Morr2}.
	
	\vspace{.1cm}
	
	Sea $C$ un campo algebraicamente cerrado que contiene a $\bb Q_p$ y es completo respecto a un valor absoluto no arquimediano $|\cdot|_C:C\to \bb R_{\geq0}$ que extiende el valor absoluto $p$-\'adico en $\bb Q_p.$ (Un ejemplo es $\bb C_p,$ los complejos $p$-\'adicos).
	
	\begin{dfn}
		El \textgravedbl tilt \textacutedbl $F=C^\flat$ de $C$ es un campo algebraicamente cerrado que contiene a $\bb F_p$ y es completo respecto a un valor absoluto no arquimediano no trivial $|\cdot|_F:F\to \bb R_{\geq0}.$
	\end{dfn}

	Como conjunto,
	$$C^\flat =\{ (a_0,a_1,...) \hspace{.2cm}| \hspace{.2cm} a_i\in C \hspace{.2cm} \text{y} \hspace{.2cm} a_i^p=a_{i-1} \}$$
	La multiplicaci\'on se define t\'ermino a t\'ermino y la suma como
	$ (a_0,a_1,...)+(b_0,b_1,...)=(c_0,c_1,...),$ donde
	$$c_i:= \lim_{i\leq n\to\infty} (a_n+b_n)^{p^{n-i}}.$$
	
	Rec\'iprocamente, sea $F$ con los hip\'otesis de la definici\'on de tilt.
	\begin{dfn}
		Un \textgravedbl untilt\textacutedbl de $F$ es un par $(C,i),$ donde $C$ es  un campo algebraicamente cerrado que contiene a $\bb Q_p$ y es completo respecto a un valor absoluto no arquimediano $|\cdot|_C:C\to \bb R_{\geq0}$ que extiende el valor absoluto $p$-\'adico en $\bb Q_p$ y $i:F\to  C^\flat$ es un isomorfismo de campos valuados.
	\end{dfn}
	
	Decimos que dos untilts $(C,i),(C',i')$ son equivalentes si existe un isomorfismo $C\cong C'$ tal que el isomorfismo inducido entre sus tilts es compatible con $i$ y $i'$, es decir, si existe $\psi: C\to C'$ isomorfismo tal que el diagrama conmuta:
	$$\xymatrix{
		C^\flat \ar^{\psi^{\flat}} [r] & C'^{\flat}  \\
		F  \ar^{i} [u]  \ar_{i'} [ur] & 
	},$$
	donde $\psi^{\flat}(a_0,a_1,...)=(\psi(a_0),\psi(a_1),...).$
	
	\begin{dfn}
		Sea $|Y_F|$ el conjunto de clases de equivalencia de untilts de $F.$
	\end{dfn}
	
	Dado un untilt $(C,i)$ de $F$ podemos construir nuevos untilts $(C,i\circ \varphi^m)$ para todo $m\in \bb Z,$ donde $\varphi$ es el automorfismo de Fronbenius.
	
	\begin{dfn}
		Decimos que dos untilts $(C,i),$ $(C',i')$ son Frobenius equivalentes si existe $m\in\bb Z$ tal que $(C,i)$ y $(C', i'\circ\varphi^m)$ son equivalentes.
	\end{dfn}
	
	Al conjunto de clases de equivalencia de untilts Frobenius equivalentes est\'a dado por el cociente $$|Y_F|/\varphi^{\bb Z},$$
	donde el grupo c\'iclico infinito $\varphi^{\bb Z}$ act\'ua en $|Y_F|$ via $$\varphi^m \cdot (C,i) = (C,i\circ \varphi^m).$$
	
	El siguiente resultado cuya demostraci\'on puede ser consultado en \cite{Morr2} nos d\'a la relaci\'on entre este conjunto de clase de equivalencia y las curva de Fargues-Fontaine.
	
	\begin{teo}
		Existe una curva $X_F^{FF}$ cuyos puntos est\'an en biyecci\'on con $|Y_F|/\varphi^{\bb Z}.$
	\end{teo}
	
	\section{La curva y las representaciones de Galois $p$-\'adicas.}
	
	En esta \'ultima secci\'on se muestra la relaci\'on de la curva de Fargues-Fontaine, espec\'ificamente sus fibrados vectoriales, con las representaciones de Galois $p$-\'adicas, haciendo uso del teorema de Harder-Narasimhan. En los fibrados vectoriales se tiene lo siguiente:

	\begin{dfn}
		Sea $X$ una superficie de Riemann. Para cualquier fibrado vectorial $E$ de $X$ asociamos dos invariantes:
		\begin{itemize}
			\item Su rango, $rk(E)\in\bb N.$
			\item Su grado, $\deg(E)\in \bb Z,$ definido como el grado de su \textgravedbl determinant line bundle\textacutedbl. El grado de un fibrado lineal $L$ se define identificando a $L$ con un divisor de Weil $\sum_{x\in X} n_x[x]$ y definimos $\deg(L)=\sum_{x\in X} n_x.$
		\end{itemize}
		De los dos invariantes anteriores se define un tercero, su pendiente 
		$$\mu(E)= \dfrac{\deg(E)}{rk(E)}\in\bb Q.$$
		
		Se dice que $E$ es semiestable si $\mu(E')\leq \mu(E)$ para todo $E'$ subfibrado vectorial de $E.$
	\end{dfn}
	
	Enunciaremos ahora el teorema principal de este secci\'on cuya demostraci\'on puede ser consultada en \cite{Morr2}.
	
	\begin{teo}[Harder-Narasimhan (H-N)]
		Sea $E$ un fibrado vectorial de $X.$ Entonces $E$ tiene una \'unica filtraci\'on de sub-fibrados
		$$0=E_0\subseteq E_1\subseteq \cdots\subseteq E_m=E,$$
		tales que:
		\begin{itemize}
			\item El fibrado cociente $\dfrac{E_i}{E_{i-1}}$ es semiestable para todo $i=1,...,m.$
			\item $\mu(\frac{E_1}{E_0})>\cdots>\mu(\frac{E_m}{E_{m-1}}).$
		\end{itemize}
	\end{teo}
	
	En general hay otras categor\'ias en las cuales existe un an\'alogo al teorema de H-N. Mencionaremos alguna de estas.
	
	\begin{itemize}
		\item Fibrados vectoriales de una curva $X.$ \\
		Sea $X$ una superficie de Riemann. Sea $Vect(X)$ la categor\'ia de fibrados en $X.$ Equipado con las nociones definidas anteriormente de rango y grado, dicha categor\'ia satisface el teorema H-N.

		\item Los pares $\bb P^1$-algebraicos completos. Antes de definir estos objetos necesitaremos la definici\'on de una funci\'on casi euclidiana.
		\begin{dfn}
			Una funci\'on casi euclidiana de un dominio entero $B$ es una funci\'on $\nu: B\to \bb N\cup \{-\infty\}$ que cumple las siguientes propiedades:
			\begin{enumerate}
				\item $\nu(f)=-\infty \iff f=0.$
				
				\item Para $f,g\in B-\{0\}$ se tiene que $\nu(f)\leq\nu(fg).$
				
				\item Si $\nu(f)=0,$ entonces $f$ es unidad.
				
				\item Si $f,g\in B$ con $\nu(g)\geq1,$ entonces existen $q,r\in B$ tales que $f=gq+r$ y $\nu(r)\leq \nu(g).$
			\end{enumerate}
		\end{dfn}  
		
		\begin{dfn}
			Un par $\bb P^1$-algebraico es un par $(B,\nu)$ que consiste  es un anillo de ideales principales $B$ y una valuaci\'on $\nu:Frac(B)\to \bb Z\cup \{\infty\} $ tal que $-\nu$ es una funci\'on casi euclidiana en $B.$ Decimos que el par el completo si 
			$$\nu (f) + \sum_{\fk p \subseteq B} ord_{\fk p}(f)=0,$$
			para todo $f\in B,$ donde $\fk p$ corre sobre todos los ideales primos distintos del cero de $B$ y $ord_{\fk p}$ denota la valuaci\'on $\fk p$--\'adica asociada en $B.$
		\end{dfn}

		Sea $(B,\nu)$ un par $\bb P^1$--algebraico completo. Un \textgravedbl fibrado vectorial en $(B,\nu)$\textacutedbl se define como un par $(M,M_\infty),$ donde $M$ es un $B$--m\'odulo libre de rango finito y $M_\infty$ es un $\ml O_\nu$--ret\'iculo dentro del espacio vectorial de dimensi\'on finita $M\otimes_B k_\nu,$ donde $k_\nu$ es la completaci\'on de $k=Frac(B)$ respecto a $\nu$  y $\ml O_\nu$ es el anillo de enteros.
		\begin{itemize}
			\item  El rango de $(M,M_\infty)$ es el rango del m\'odulo $M.$
			\item El grado se define comparando bases de $M$ y $M_\infty.$ 
		\end{itemize}
		As\'i, la categor\'ia de fibrados vectoriales de un par $\bb P^1$-algebraico $(B,\nu)$ cumple el teorema H-N. \\
		En particular, para el par $\bb P^1$-algebraico completo $(B_e,\nu_{dR})$ de la curva de Fargues-Fontaine, los $(B,\nu)$--pares solo son llamados $B$-pares. Aqu\'i un fibrado vectorial es un par $(M,M_{dR})$ tal que
		\begin{itemize}
			\item $M$ es un $B_e$--m\'odulo libre de rango finito.
			\item $M_{dR}$ es un $B_{dR}^+$--ret\'iculo dentro de $M\otimes_{B_e} B_{dR}.$
		\end{itemize}
		La siguiente proposici\'on relaciona las dos categor\'ias antes mencionadas. Su demostraci\'on puede ser consultada en \cite{Morr2}.
		
		\begin{prop}
			La categor\'ia de $(B_e,\nu_{dR})$--pares se identifica con la categor\'ia $Vect(X^{FF})$ de fibrados vectoriales de la curva de Fargues--Fontaine.
		\end{prop}
		
		\item  Espacios vectoriales con filtraciones.\\
		Dada una extensi\'on de campos $L/F,$ sea $VectFil_{L/F}$ la categor\'ia de pares $(V, Fil^{\bullet}V_L),$ donde $V$ es un $F$--espacio vectorial de dimensi\'on finita y $Fil^{\bullet}$ es una filtraci\'on en $V_L=V\otimes_F L$ separada y exhaustiva.
		\begin{itemize}
			\item $rk(V,Fil^{\bullet}V_L)=\dim_F V.$
			\item $\deg(V,Fil^{\bullet}V_L)=\sum_{i\in \bb Z} i\dim_L(gr^i V_L)$
		\end{itemize}
		As\'i, esta categor\'ia cumple con el teorema H-N.
		
		\item Isocristales.\\
		Sea $K$ campo perfecto de caracter\'istica $p$ y $K_0=Frac(W(K)).$ Un isocristal sobre $K$ es un par $(D,\varphi_D),$ donde $D$ es un $K_0$--espacio vectorial de dimensi\'on finita y $\varphi_D:D\to D$ es un isomorfismo $\varphi$--semilineal (esto es, $\varphi_D(ad)=\varphi(a)\varphi_D(d)$ para todo $a\in K_0$ y $d\in D.$)
		\begin{itemize}
			\item $rk(D,\varphi_D)=\dim_{K_0}D.$
			\item $\deg(D,\varphi_D)=-\deg^+\det(D,\varphi_D),$ donde $\deg^+$ de un isocristal de rango uno $(L,\varphi_L)$ se define escogiendo un elemento b\'asico $e\in L,$ luego $\varphi_L(e)=ae$ para alg\'un $a\in K_0$ y $\det^+(L,\varphi_L):= \nu_p(a),$ ya que $K_0$ es campo de valuaci\'on discreta.
		\end{itemize}
		As\'i, esa categor\'ia cumple con el teorema de H-N. A esta categor\'ia se le suele denotar como $\varphi-Mod_{K_0}.$\\
		
		Sea $\lambda=\dfrac{d}{h}\in\bb Q,$ donde $(d,h)=1$ y $h>0.$ Podemos definir un isocristal $(D_\lambda,\varphi_\lambda)\in \varphi-Mod_{K_0}$ como
		\begin{itemize}
			\item $D_\lambda=K_0^h,$ con elementos b\'asicos $e_1,...,e_h.$
			\item $\varphi_\lambda:D_\lambda\to D_\lambda$ como el \'unico endomorfismo semilineal que satisface:
			$$\varphi_\lambda(e_i)= \left\{ \begin{array}{lcc}
				e_{i+1}, & si & i=1,...,h-1, \\
				\\ p^{-d}e_1, & si  &  i=h.
			\end{array}
			\right.$$
		\end{itemize}
		El isocristal $(D_\lambda,\varphi_\lambda)$ tiene rango $h,$ grado $d$ y pendiente $\lambda.$
		
		\begin{lem}
			Sea $(D,\varphi_D)\in\varphi-Mod_{\bb Q_p}$ isocristal sobre $\bb F_p.$ Entonces el par 
			$$((B_{cris}\otimes_{\bb Q_p} D)^{\varphi=1}, B^+_{dR}\otimes_{\bb Q_p} D)$$
			es un fibrado vectorial $(B_e,\nu_{dR})$ con rango y grado dado por el rango y grado del isocristal. 
		\end{lem}

		Este lema asocia a un isocristal fibrado un  fibrado vectorial en $X^{FF},$ denotado por $\xi(D,\varphi_D).$ En particular, en el caso del isocristal $(D_\lambda,\varphi_\lambda)$ con $\lambda\in \bb Q,$ escribimos $\ml O_{X^{FF}}(\lambda):=\xi(D_\lambda,\varphi_\lambda).$ Por el lema, $\ml O_{X^{FF}}(\lambda)$ tiene rango $h,$ grado $d$ y pendiente $\lambda$ si $\lambda=\dfrac{d}{h}$ con $(d,h)=1$ y $h>0.$
		\begin{teo}[De clasificaci\'on de fibrados vectoriales en $X^{FF}.$]
			Sea $E$ un fibrado vectorial en $X^{FF}.$ Entonces existe una \'unica sucesi\'on de racionales $\lambda_1\geq \cdots\geq \lambda_m$ tales que $E$ es isomorfo a 
			$$\bigoplus_{i=1}^m \ml O_{X^{FF}}(\lambda_i) =\bigoplus_{i=1}^m \xi(D_{\lambda_i},\varphi_{\lambda_i}).$$
		\end{teo}
		
		\begin{cor}\hspace{.2cm}
			\begin{itemize}
				\item El funtor $\xi(-):\varphi-Mod_{K_0}\to Vect(X^{FF})$ es esencialmente sobreyectivo.
				
				\item Sea $E$ un fibrado vectorial de $X^{FF}$ y $\lambda\in\bb Q.$ Entonces $E$ es semiestable de pendiente $\lambda$ si y s\'olo si $E$ es isomorfo a $\ml O_{X^{FF}}(\lambda)^m$ para alg\'un $m\geq1.$
				
				\item La categor\'ia de fibrados vectoriales de $X^{FF}$ semiestables y de pendiente cero  equivalente a la categor\'ia de $\bb Q_p$--espacios vectoriales de dimensi\'on finita v\'ia $V\mapsto V\otimes_{\bb Q_p} \ml O_{X^{FF}}.$
			\end{itemize}
		\end{cor}

		\item $\varphi-ModFil_{K/K_0}.$\\
		Esta es la categor\'ia de los tripletes $(D,\varphi_D,Fil^{\bullet}D_K)$ donde
		\begin{itemize}
			\item $D$ es un $K_0$--espacio vectorial de dimensi\'on finita.
			
			\item $\varphi_D:D\to D,$ isomorfismo $\varphi$-lineal.
			
			\item $Fil^{\bullet}$ es una filtraci\'on en $D_K:=D\otimes_{K_0} K$ separada y exhaustiva.
		\end{itemize}
		Es decir, $(D,\varphi)\in \varphi-Mod_{K_0}$ y $(D,Fil^{\bullet}D_K)\in VectFil_{K_0/K}.$ Definimos:
		\begin{itemize}
			\item $rk(D,\varphi_D,Fil^{\bullet}D_K):=\dim_{K_0} D.$
			\item $\deg(D,\varphi_D,Fil^{\bullet}D_K):=\deg(D,\varphi_D) + \deg(D,Fil^{\bullet} D_K).$
		\end{itemize}
		De esta forma, esta categor\'ia cumple con el teorema de H-N.
	\end{itemize}
	
	Como $B_{cris}\subseteq B_{dR},$ si $V$ es representaci\'on cristalina entonces $V$ es de deRham y $D_{B_{dR}}(V)= D_{B_{cris}}(V)\otimes_{K_0} K.$ As\'i el funtor 
	$$D_{cris}:Rep_{cris}(G_K)\to\varphi-ModFil_{K/K_0}$$
	es completamente fiel.
	
	Fontaine conjetur\'o que todo isocristal filtrado estaba en la imagen esencial de $Rep_{cris}(G_K)$ si y s\'olo si era d\'ebilmente admisible \cite[Cap\'itulo 8: classification des fibr\'{e}s vectoriels: le cas $F$ alg\'{e}briquement clos]{Far} , que en este caso, es equivalente a pedir que sea semiestable y de pendiente cero. Esto demuestra que 
	$$D_{cris}:Rep_{cris}(G_K)\to\varphi-ModFil_{K/K_0}^{w.a.},$$
	donde $w.a.$ denota \textgravedbl weakly admissible\textacutedbl (d\'ebilmente admisible) es una equivalencia de categor\'ias. Tambi\'en es posible demostrar que la categor\'ia de representaciones semiestables se puede describir como isocristales filtrados d\'ebilmente admisibles con operador de monodrom\'ia \cite{Far}.
	
	Finalmente, hemos llegado a uno de los objetivos de la teor\'ia de Hodge $p$-\'adica, describir clases de representaciones de Galois $p$-\'adicas en t\'erminos de pura \'algebra lineal.  
	
	\section{Ep\'ilogo}
	
	El estudio de la curva de Fargues-Fontaine es actualmente un tema de investigaci\'on muy activo dentro del \'area de la geometr\'ia aritm\'etica. Trabajos como los de Peter Scholze de perfectoides \cite{scholze2014perfectoid} o de Fargues respecto al programa de Langlands \cite{fargues2021geometrization} son muestra del actual avance en este t\'opico. Otro punto de vista que particularmente estamos estudiando son los sistemas din\'amicos de funciones racionales sobre la curva: estudiar puntos peri\'odicos o los conjuntos de Fatou y Julia podr\'ian llevarnos a descubrir a\'un m\'as informaci\'on acerca de este objeto matem\'atico tan importante.

	\subsubsection*{{\bf\em\scriptsize AGRADECIMIENTOS}} El primer autor expresa su gratitud a CONAHCyT por la beca obtenida para estudiar el doctorado en el CIMAT y con el cual fue posible el desarrollo del presente art\'iculo.

	\newpage
	
	\begin{tabular}[t]{l}
		\emph{Direcci\'on del autor:}\\
		{Centro de Investigaci\'on en Matem\'aticas (CIMAT),}\\
		{Unidad M\'erida,}\\
		Parque Cient\'ifico y Tecnol\'ogico de Yucat\'an
		Km 5.5 Carretera Sierra Papacal \\
		{Chuburn\'a Puerto Sierra Papacal; CP 97302, M\'erida, Yucat\'an.}\\
		e-mail: \texttt{jorge.robles@cimat.mx}
	\end{tabular}

     \begin{tabular}[t]{l}
     	\emph{Direcci\'on del coautor:}\\
     	{Centro de Investigaci\'on en Matem\'aticas (CIMAT),}\\
     	{Unidad M\'erida,}\\
     	Parque Cient\'ifico y Tecnol\'ogico de Yucat\'an
     	Km 5.5 Carretera Sierra Papacal \\
     	{Chuburn\'a Puerto Sierra Papacal; CP 97302, M\'erida, Yucat\'an.}\\
     	e-mail: \texttt{rogelio.perez@cimat.mx}
     \end{tabular}

	\label{artuno}

\end{document}